\magnification=1100
\baselineskip=16pt

\hsize=150mm
\vsize=205mm

\let \al=\alpha
\let \be=\beta
\let \var=\varphi
\let \vare=\varepsilon

\let \th=\theta
\let \la=\lambda

\let \q=\quad

\let \med=\medskip
\let \smal=\smallskip

\def\R{{\rm I\kern
-1.6pt{\rm R}}}
\def\C{{\rm |\kern
-4.6pt{\rm C}}}
\def\N{{\rm I\kern
-4.0pt{\rm N}}}

\def\eq{\eqno}
\def\q{\quad}

\def\ter{\hfill \vrule width 5 pt height 7 pt depth - 2 pt\smallskip}

\

{\bf \centerline{ Persistence, permanence and global stability} 

\smal
\centerline{for an $n$-dimensional Nicholson system }}

\

\vskip0.3in

\centerline{Teresa Faria\parindent=0cm{\footnote {$^{a}$}{Departamento de Matem\'atica  and CMAF, Faculdade de
Ci\^encias, Universidade de Lisboa, Campo Grande, 1749-016 Lisboa,
Portugal}}
\parindent=0cm{\footnote{$^{\star}$}{Corresponding author.
Fax:~+351 21 795 4288, tel:~+351
21 790 4929, e-mail:~tfaria@ptmat.fc.ul.pt.}}
and Gergely R\"ost\parindent=0cm{\footnote{$^b$}{Bolyai Institute, University of Szeged,
 Aradi v\' ertan\' uk tere 1., H-6720 Szeged, Hungary}}
 }

\vskip0.3in

{\it Suggested running head: Global asymptotic behavior for Nicholson's systems}

\vskip0.3in

\centerline{\bf Abstract}
\vskip0.3in

For a Nicholson's blowflies system with patch structure and multiple discrete delays, we analyze several features of the global asymptotic behavior of its solutions. It is shown that if the spectral bound of the community matrix is non-positive, then the population becomes extinct on each patch, whereas the total population uniformly persists if the spectral bound is positive. 
Explicit uniform lower and upper bounds for the asymptotic behavior of solutions are also given.
When the population uniformly persists,  the existence of a unique positive equilibrium is established, as well as a sharp criterion for its absolute global asymptotic stability, improving  results in the recent literature. While our system is not cooperative, several sharp threshold-type results about its dynamics are proven, even when the community matrix is reducible, a case usually not treated in the literature.

\vskip0.3in


\vskip0.3in


\vfill\eject

{\bf 1. Introduction}

\

In recent years, population dynamics  models with patch structure and delays have attracted the attention of an increasing number of mathematicians and biologists. 
The heterogeneity of the  
environment is inherently captured by patchy models, in which the  
spatial distribution of the population is governed by both the  
migration between patches and the growth of the local populations,  
which depends on the resources of each particular patch.  
Patch-structured systems of differential equations are also used as  
disease models with transitions between stages of normal and infected  
cells. Delay differential equations (DDEs)  frequently  
provide quite realistic  models in population dynamics, epidemiology  
and mathematical biology in general, since the  incorporation of  
delays appears naturally to express the maturation period of  
biological species, the maturation time of blood cells,  the  
incubation period in disease models,  and  several other features. Understanding the interplay of spatial dispersal and time delays is  therefore a key point for many models.

\med

In the present paper, we study some aspects of the asymptotic behavior of solutions for the following Nicholson's blowflies system with patch structure and multiple discrete delays:
$$x_i'(t)=-d_ix_i(t)+\sum_{j=1}^n a_{ij}x_j(t)+\sum_{k=1}^m \be_{ik}    x_i(t-\tau_{ik})e^{-x_i(t-\tau_{ik})},\q i=1,\dots,n,\eq(1.1)$$
where $d_i>0, a_{ij}\ge 0,\tau_{ik}> 0,  \be_{ik} \ge 0$ and
$$\be_i:=\sum_{k=1}^m \be_{ik} >0\eq(1.2)$$
for $i,j=1,\dots, n, k=1,\dots,m$.
By condition (1.2),  there is at least one delayed nonlinearity on each patch $i$.
  To simplify the notation and without loss of generality, in what follows we shall always assume that $a_{ii}=0$ for all $1\le i\le n$.

 Among other applications, system (1.1) fits as a population model for the growth of  single or multiple biological species
divided into  $n$ patches or classes, with migration of the populations among them. 
On each patch $i$, $x_i(t)$ denotes the density of the population, $d_i$ is its decreasing rate,  
the birth function  is of Nicholson-type 
$\sum_{k=1}^m \be_{ik}    x_i(t-\tau_{ik})e^{-x_i(t-\tau_{ik})}$, and the coefficients $a_{ij}$ are the migration rates of populations moving from patch $j$ to patch $i$.
 In view of this biological meaning, it is natural to take
$$d_i=m_i+\sum_{j=1}^n a_{ji},\q m_i>0,
\eq(1.3)$$
where  $m_i$ is the mortality rate on patch $i$. Therefore, together with conditions $a_{ii}=0$ and (1.2),  unless otherwise stated,  in what follows we  assume  (1.3).

Model (1.1) was motivated by the  celebrated scalar Nicholson's blowflies  equation 
$$x'(t)=-dx(t)+\be x(t-\tau)e^{-a x(t-\tau)},$$
where $d,\be , a,\tau >0$,
 introduced by Gurney et al.~[6] in 1980 as a model for the Australian sheep-blowfly population,
as it agreed  with the Nicholson's experimental data published in [13]. Since then, Nicholson's equation has been generalized, modified, and extensively  studied by many mathematicians, in what concerns  stability, persistence, existence and attractivity of periodic or almost periodic solutions, occurrence of bifurcations, and other dynamical aspects.
In contrast, the literature on Nicholson's systems is quite recent and scarce. We refer to the works of  Liu [10,~11],  Berezansky et al.~[1], Faria [3], Liu and Meng [12], and  Wang [18].

 Throughout the paper, we  designate $A,B,D$ as the matrices
$$A=[a_{ij}],\q D=diag\, (d_1,\dots,d_n),\q B=diag\, (\be_1,\dots,\be_n),
\eq(1.4)$$
and refer to
$$M:=A+B-D$$ as 
the {\bf community matrix}. 
  The algebraic properties of the community matrix will play an important role in the study of either the persistence or the extinction of the species in all  patches, as well as in the existence  of a positive equilibrium -- whereas the stability of the positive equilibrium  depends heavily on the shape of the non-linear terms in (1.1). While most  papers dealing with multiple dimensional DDEs used in population dynamics only consider the situation of an {\it irreducible} community matrix,  in the present paper we also treat the case a {\it reducible}  matrix. 

The present paper is as an extension of the research in [3], where sufficient conditions for the global attractivity of both the trivial equilibrium and the positive equilibrium, when it exists,  were established.  Here, we pursue a deeper analysis of system (1.1), improving the criteria established in [3]  and addressing new aspects of its dynamics. The paper provides answers for current important open problems. Namely, it gives a threshold condition  for the extinction of the populations in all patches versus the uniform persistence of the total population -- which applies even for the particular case of a {\it reducible} community matrix --,  shows the existence of a positive equilibrium under very general assumptions, and establishes a (sharp) criterion for its absolute global asymptotic stability. Some of our results
naturally hold  for delayed systems  with a more general class of nonlinearities, however  the criteria for the global asymptotic stability of the positive equilibrium, as well as 
some explicit upper and lower bounds for the asymptotic behavior of solutions 
are very specific to the Ricker-type nonlinearity in (1.1).

Some of main techniques used here rely on  M-matrix theory and on properties of cooperative systems of DDEs. We refer the reader to the 
  monograph  of Fiedler [5] for properties of M-matrices, 
the 
  monograph  of Smith on monotone systems [15]  for   cooperative behavior of DDEs, and the recent book of Smith and  Thieme [16] for terminology and results on population persistence. Also, the  method developed by Faria and Oliveira [4]  to study the stability of linear $n$-dimensional DDEs was used   to address the local asymptotic stability of the equilibria of system (1.1), an aspect previously exploited in [3]. Another major   source of inspiration for our work was  
the paper of Hofbauer [8], where   the concept of {\it saturated equilibrium} for autonomous systems of ordinary differential equations (ODEs) which are positively invariant in the positive cone $\R^n_+$ was introduced, and  powerful results on the existence of a saturated equilibrium for dissipative systems were established. Hofbauer's results were a key point in our research, to provide  a very general criterion for the existence of a unique positive fixed point of (1.1).

\smal
We now introduce some notation and set some terminology.
For the DDE (1.1),  we choose the usual phase space $C:=C([-\tau,0];\R^n)$ of continuous functions from $[-\tau,0]$ to $\R^n$ with the supremum norm $\|\var\|=\max_{\th\in [-\tau,0]}|\var(\th)|$, where $\tau=\max_{1\le i\le n,1\le k\le m} \tau_{ik} >0$ and $|\cdot|$ is any chosen norm in $\R^n$. In Section 2, when dealing with the concept of {\it $\rho$-uniform persistence}, for practical reasons it will be convenient to choose the norm $|x|=\sum_{i=1}^{n}|x_{i}|$, for the calculations
in the persistence proof. For similar reasons, in Section 5 we choose the maximum norm in $\R^n$, to address the global asymptotic stability of the positive equilibrium.
Due to the biological interpretation of model (1.1), we shall restrict our attention to non-negative solutions, and consider as  set of admissible initial conditions  either  the positive cone $C^+=\{ \var\in C: \var_i(\th)\ge 0$ for all $\th\in[-\tau,0],i=1,\dots,n\}$ or the  subset $C_0^+$ of $C^+$ of  functions which are strictly positive at zero, $C^+_0=\{ \var\in C^+:\var_i(0)>0, i=1,\dots,n\}$. One can use the method of steps to verify that both sets  $C^+$ and
  $C_0^+$ are positively invariant under (1.1). Moreover,
 for each $\var\in C^+$ system (1.1) has a unique solution $x(t)=x(t;\var)$ defined on $[0,\infty)$,
 with $x_i(t)$   positive on $[0,\infty)$ provided that $x_i(0)=\var_i(0)>0$. As usual, segments of solutions in the phase space $C$ are denoted by $x_t$, $x_t(\th)=x(t+\th), \th\in[-\tau,0]$, with components $x_{t,i}$.
 When analyzing (1.1), our concept of  stability always refers to the setting of {\it admissible solutions}, i.e., solutions  $x(t;\var)$ with $\var$ in the set of admissible initial conditions. In particular, the trivial equilibrium of (1.1) is {\it globally asymptotically stable} (GAS) if it is  stable and  attracts  all  solutions $x(t)=x(t;\var)$ of (1.1)  with initial conditions $\var\in C^+$, i.e., $\lim_{t\to \infty} x(t)=0$;
 if  $x^*>0$ is an equilibrium of (1.1), $x^*$ is said to be GAS if it is  stable and  attracts  all  solutions $x(t)=x(t;\var)$ of (1.1)  with initial conditions $\var\in C^+_0$.

For a vector $c\in\R^n$, we also use $c$ to denote the constant function  $\var (\th)=c$ for $\th\in [-\tau,0]$ in $C$.
 A  vector $c$ is said to be {\it positive}, or {\it non-negative}, if all its components are positive, or non-negative, respectively.  We define in a similar way {\it positive} and {\it non-negative functions} in $C$, and  {\it positive} and {\it non-negative matrices}.

  \smal
 
 We  recall  below some  concepts  from matrix theory, included here for convenience of the reader, since  they will be often referred to  in the next sections.

 \smal
 
 {\bf Definition 1.1}. Let $N=[n_{ij}]$ be an ${n\times n}$ matrix. We say that 
   $N$ is {\bf cooperative} if  its  off-diagonal entries are non-negative: $n_{ij}\ge 0$ for $j\ne i$.
The matrix
  $N$ is a {\bf reducible matrix} if there is a simultaneous permutation of rows and columns that brings $N$ to the form
  $$\left [\matrix{ N_{11}&0\cr N_{21}&N_{22}\cr}\right]$$
  with $N_{11}$ and $N_{22}$   square matrices;
     $N$ is an {\bf irreducible matrix} if it is not reducible.
 The {\bf spectrum} of $N$ is denoted by $\sigma (N)$. The
{\bf spectral bound}  of $N$ is defined as  
$$s(N)=\max\{ Re\, \la: \la  \in\sigma (N)\}.$$
The matrix
  $N$ is said to be an {\bf M-matrix} if $a_{ij}\le 0$ for $i\ne j$ and all its eigenvalues have non-negative real parts. If $N$ is an M-matrix and $\det N\ne 0$, then we say that $N$ is a {\bf non-singular M-matrix}.
  
  \smal

  It is well-known that there are several equivalent ways of defining M-matrices and non-singular M-matrices, see e.g. [5, 17] for further properties of these matrices. However we emphasize that many authors use the term {\it M-matrix} with the above meaning of the term  {\it non-singular M-matrix}.
We also recall that if a square matrix $N$  is cooperative and irreducible, then its spectral bound $s(N)$ is always a simple, dominant eigenvalue, with a positive associated eigenvector [17].

\smal

The remainder  of the paper consists of four sections. The  persistence and permanence of the Nicholson-type system (1.1), two crucial aspects in population dynamics (see e.g. [16]), are studied in Section 2.   When $s(M)>0$, a further  analysis is carried out to obtain strong uniform persistence of the population at least on one patch, and for all the patches in the case of  an irreducible community matrix.
Explicit lower and upper uniform bounds for the positive solutions of (1.1) given in terms of the coefficients in (1.1) are also included.
In Section 3, we prove the global attractivity of the equilibrium 0 when $s(M)\le 0$, which means the extinction of the populations in all patches. Therefore, a threshold criterion for extinction versus persistence  is provided; moreover, this persistence is uniform in
  the special case of an irreducible community matrix. Clearly, from the point of view of applications, it is most relevant to study the existence, stability and attractivity of a positive equilibrium. The last sections are dedicated to these aspects.  In Section 4, we study the undelayed ODE version of (1.1), obtained  by taking all the delays equal to zero, and prove the existence of a unique positive equilibrium for (1.1) if $Mc>0$ for some positive vector $c$.
  Finally, in Section 5 we give a sharp criterion for the absolute global asymptotic stability of such equilibrium, which   significantly improves recent results in the literature, see e.g. [1,~3,~10,~11].

 \
 
 {\bf 2. Boundedness of solutions, persistence, permanence}
 
 \med
 
In this section, we analyze the permanence and persistence of (1.1).

 We first observe that
  condition (1.3)  implies that  the matrix $D-A^T$ is diagonally dominant, therefore from Theorems 5.14 and 5.1 in  [5] it follows that $D-A^T$ is always a non-singular M-matrix, and thus  $D-A$ as well. As an immediate consequence of $D-A$ being a non-singular M-matrix, we get 
 the boundedness of all admissible solutions of (1.1). 

\proclaim{Theorem 2.1}. System (1.1) is dissipative on $C^+$, i.e., the components of all solutions of (1.1) with initial conditions in $C^+$ are uniformly bounded. To be more precise, 
all the solutions $x(t)=x(t,\var)$ of (1.1) with initial conditions $x_0=\var\in C^+$ satisfy
$$d_iu_i- \sum_{j=1}^n a_{ij}\ u_j\le \be_i e^{-1},\q i=1,\dots,n,
\eq(2.1)$$
or, in other words,
$$\left [\matrix{u_1\cr\vdots\cr u_n\cr}\right]\le (D-A)^{-1}\left [\matrix{\be_1\cr\vdots\cr\be_n\cr}\right ]e^{-1},
\eq(2.2)$$
where $u_i=\limsup_{t\to \infty}x_i(t),\ i=1,\dots,n$.

{\it Proof}.  Fix $s>0$. For any $\var\in C^+$, consider the solution  $x(t)=x(t,\var)$ of (1.1), and define 
$\bar u_i=\sup_{t\in [0,s]} x_i(t),\ i=1,\dots,n$. Since  $h(x):=xe^{-x}\le e^{-1}, x\ge 0$, then
$x_i'(t)\le -d_ix_i(t)+\sum_{j=1}^n a_{ij}\bar u_j+\be_i e^{-1}$, implying that
$e^{d_it}x_i(t)\le x_{0i}+(e^{d_it}-1)\bar \eta_i /d_i, 0\le t\le s$, where $\var(0)=(x_{01},\dots,x_{0n})\in\R^n_+$, and $\bar\eta_i=\sum_{j=1}^n a_{ij}\bar u_j+\be_i e^{-1}$.  Hence
we obtain 
$$x_i(t)\le x_{0i}e^{-d_it}+d_i^{-1}\bar\eta_i (1-e^{-d_it}),\q i=1,\dots,n,
\eq (2.3)$$
 from which we deduce $d_i\bar u_i\le d_ix_{0i}+\be_ie^{-1}+\sum_{j=1}^n a_{ij}\bar u_j,\, i=1,\dots,n$; in other words, for $\bar u=(\bar u_1,\dots,\bar u_n)$, we have
$$(D-A)\bar u\le c,\q {\rm with}\q c=\left [\matrix{d_1x_{01}\cr\vdots\cr d_nx_{0n}\cr}\right ]+\left [\matrix{\be_1\cr\vdots\cr\be_n\cr}\right ]e^{-1}.
\eq(2.4)$$
Since $D-A$ is a non-singular M-matrix, then  its inverse is a non-negative matrix [5,~Theorem 5.1], and from (2.4) we get
$\bar u\le (D-A)^{-1}c.
$
This estimate does not depend on $s>0$, thus we derive 
$$u\le (D-A)^{-1}c,
\eq(2.5)$$
for $u=(u_1,\dots,u_n)$ and $u_i=\limsup_{t\to \infty}x_i(t),\ i=1,\dots,n$, implying that all positive solutions are bounded. Next, we prove that the uniform estimate (2.2) holds.

Let $\vare>0$. For $t>0$ large, we have $x_i(t)\le u_i+\vare$, thus the estimate (2.3) is obtained with $\bar\eta_i$ replaced by
$\eta_i=\sum_{j=1}^n a_{ij}\ (u_j+\vare)+\be_i e^{-1},$ for $i=1,\dots,n$.
By letting $\vare\to 0^+$ and $t\to\infty$, it follows that $d_iu_i\le \be_i e^{-1}+\sum_{j=1}^n a_{ij}\ u_j$, for all $i$, which proves (2.1), and therefore 
$(D-A)u\le [\be_1\ \cdots\ \be_n]^Te^{-1}$.\ter

\med

For the definitions of persistence and permanence given below, see e.g.   [9].

\med

{\bf Definition 2.1.} System (1.1) is said to be {\bf persistent} (in $C^+_0$) if any  solution $x(t;\var)$ with initial condition $\var\in C^+_0$ is bounded away from zero, i.e.,
$\liminf_{t\to\infty} x_i(t;\var)>0, 1\le i\le n$, for any any $\var\in C^+_0$; and {\bf uniformly persistent} (in $C^+_0$) if there is $\eta >0$ such that $\liminf_{t\to\infty} x_i(t;\var)\ge \eta, 1\le i\le n$, for any any $\var\in C^+_0$.
System (1.1) is said to be {\bf permanent} (in $C^+_0$) if there are positive constants $m_0,M_0,$ with $m_0<M_0,$ such that, given any $\var\in C^+_0$, there exists $t_0=t_0(\var)$  such that $m_0\le x_i(t,\var)\le M_0$ for $1\le i\le n$ and $t\ge t_0$. 

\med

The notion of persistence in Definition 2.1 means that the population persistence on each patch. In the following, we shall discuss population persistence on a particular patch, on a given subset of patches, or the persistence of the total population. In order to perform such analysis, we also use the more general terminology of $\rho$-persistence as it has been presented in the monograph of Smith and Thieme [16].

\med

{\bf Definition 2.2.} Let $X$ be a nonempty set of a Banach space and $\rho:X\to\R_{+}$.
A semiflow $\Phi:\R_{+}\times X\to X$  is called {\bf uniformly
weakly $\rho$-persistent}, if there exists some $\varepsilon>0$
such that 
$$
\limsup_{t\to\infty}\rho(\Phi(t,x))>\varepsilon\qquad\forall x\in X,\ \rho(x)>0.
$$
$\Phi$ is called {\bf uniformly (strongly) $\rho$-persistent}
if there exists some $\varepsilon>0$ such that 
$$
\liminf_{t\to\infty}\rho(\Phi(t,x))>\varepsilon\qquad\forall x\in X,\ \rho(x)>0.
$$

System (1.1) generates a semiflow on $C^{+}$. To discuss the persistence
on a given patch $j$, we may choose $\rho_j(\phi):=\phi_{j}(0)$. Then the uniform $\rho_j$-persistence of (1.1) for all $j$ coincides with the concept of uniform persistence of (1.1) in the sense of Definition 2.1. Choosing
$\rho(\phi):=|\phi(0)|=\sum_{i=1}^{n}\phi_{i}(0)$, we can talk about
the persistence of the total population of (1.1).

\med

Next, we prove the persistence of system (1.1).

\proclaim {Theorem 2.2}. Consider (1.1) and assume that there is  a vector $c=(c_1,\dots,c_n)>0$ such that
$$\be_ic_i>d_ic_i-\sum_{j=1}^n a_{ij}c_j,\q i=1,\dots,n.\eq(2.6)$$
Then, $\liminf_{t\to\infty} x_i(t;\var)>0, 1\le i\le n$, for any  solution $x(t;\var)$ with initial condition $\var\in C^+_0$.

{\it Proof}.  The statement was proved in [3, Lemma 2.5], with (2.6) replaced by the condition $\be_i>d_i-\sum_{j=1}^n a_{ij}$ for all $i$. The proof of this theorem is similar 
after the changes of variables $x_i\mapsto c_i^{-1}x_i, 1\le i\le n$, so it is omitted.\ter

\med

Clearly the matrix $M$ is  cooperative.
Note that condition (2.6) is equivalent to saying that $Mc>0$, for some positive vector $c$. If the matrix $A$ is irreducible, the matrix $M$ is irreducible as well, thus the spectral bound of $M$, $s(M)=\max \{ Re\, \la: \la \in\sigma(M)\}$, is an eigenvalue of $M$ with  a positive associated eigenvector, and (2.6) holds. Actually, for irreducible matrices one can use algebraic arguments -- or, in alternative, the results in Section 3  (cf. Theorem 3.3) --  to show that the converse is also true. Hence, for irreducible matrices, $s(M)>0$ is a criterion for the persistence of (1.1) in $C_0^+$, which will be shown  to be  sharp. For the  reducible case, however,  $s(M)>0$  is not a sufficient condition for persistence, as shown by the following counter-example.

\smal

{\bf Example 2.1}. Consider the 2-patch system 
$$
 \eqalign{
x_1'(t)&=-d_1 x_1(t)+\be_1e^{-x_1(t-\tau_1)}x_1(t-\tau_1)\cr
x_2'(t)&=-d_2 x_2(t)+\be_2e^{-x_1(t-\tau_2)}x_2(t-\tau_2)+a_{21}x_1(t)\cr}\eq(2.7)
$$
with $\be_1,\be_2,d_1,d_2, a_{21}>0, \tau_1,\tau_2\ge 0$, and
$\be_1<d_1$,  $\be_2>d_2$. 
Then we have
$$A=\left [
\matrix{0&0\cr
a_{21}&0\cr}
\right ], \ M=\left [\matrix{\be_1-d_1&0\cr
a_{21}&\be_2-d_2\cr}\right ],$$
 so $s(M)=\be_2-d_2>0$. On the other hand the first equation of (2.7) decouples, and since $\be_1<d_1$, we can apply Proposition 3.1 of [14] to the scalar equation of $x_1(t)$ to see that $x_1(t)\to 0$ as $t\to\infty$ for all values of the delay $\tau_1$, so (2.7) is not persistent.

\med

To study the permanence of (1.1), we start with an auxiliary lemma.

\proclaim {Lemma 2.1}. Consider the system
$$x_i'(t)=-d_ix_i(t)+\sum_{j=1}^n a_{ij}x_j(t)+\sum_{k=1}^m \be_{ik}    x_i(t-\tau_{ik})e^{-c_ix_i(t-\tau_{ik})},\q i=1,\dots,n,\eq(2.8)$$
where $c_1,\dots,c_n>0$, all the other coefficients are as in (1.1), and conditions (1.2) and (1.3) hold. Assume in addition that
$${\bf (A1)}\hskip 2cm\gamma_i:={\be_i\over {d_i-\sum_{j=1}^n a_{ij}}}>1,\q i=1,\dots,n. \hskip 2cm
\eq(2.9)$$
Let $t_*\ge 0, L>1$, and $x(t)$ be a positive solution of (2.8) satisfying $x_i(t)\le L$ for $t\ge t_*$ and $i=1,\dots,n$. Choose $m>0$ such that
$$c_im<1,\q h_i(m)\le h_i(L)\q {\rm and}\q e^{c_im}\le \gamma_i,\q i=1,\dots,n,
\eq(2.10)$$
where $h_i(x)=xe^{-c_ix}, x\ge 0,$ for $1\le i\le n$.
Then $\liminf_{t\to\infty} x_i(t)\ge m$ for all $1\le i\le n$.

{\it Proof}. The proof was inspired by an idea in [2]. Let $x(t)$ be a solution of (2.8), and  fix $m$ satisfying (2.10). Note that each function $h_i$ is strictly increasing on $[0,c_i^{-1}]$ and strictly decreasing on $[c_i^{-1},\infty)$. First, we prove:

\med

{\it Claim 1}. If $\displaystyle \min_{1\le j\le n,t\in [T,T+\tau]} x_j(t)\ge m $ for some $T\ge t_*$, then $ x_j(t)\ge m $ for all $t\ge T$ and $j=1,\dots, n$.

\med

Without loss of generality take $t_*=T=0$, and assume that $x_j(t)\ge m$ for $t\in [0,\tau]$ and $j=1,\dots, n$. Let $t_0\in [\tau, 2\tau]$ and $i\in \{1,\dots, n\}$ such that $x_i(t_0)= \min_{1\le j\le n,t\in [\tau,2\tau]} x_j(t)$.

If $x_i(t_0)<m$, we have 
$$ 0\ge x_i'(t_0)=-d_ix_i(t_0)+\sum_{j=1}^n a_{ij} x_j(t_0)+\sum_{k=1}^m \be_{ik} h_i(x_i(t_0-\tau_{ik})).
\eq (2.11)$$
Note that $x_i(t_0-\tau_{ik})\in [m,L]$ if $t_0-\tau_{ik}\in [0,\tau]$, and $x_i(t_0-\tau_{ik})\ge x_i(t_0)$ if $t_0-\tau_{ik}\in [\tau, t_0]$, hence $h_i(x_i(t_0-\tau_{ik}))\ge \min \{ h_i(x_i(t_0)), h_i(m)\} =h_i(x_i(t_0))$, and from $e^{c_im}\le \gamma_i$ we obtain
$$0\ge \left (-d_i+\sum_{j=1}^n a_{ij}+\be_ie^{-c_ix_i(t_0)}\right ) x_i(t_0) > \left (-d_i+\sum_{j=1}^n a_{ij}+\be_ie^{-c_im}\right )x_i(t_0)\ge 0,
$$
and a contradiction. Thus, $x_i(t_0)\ge m$. By the method of steps, this proves Claim 1.

\med

Next, denote $s_0:=\min_j \min_{t\in[0,\tau]} x_j(t)>0$.

If $s_0\ge m$, the result follows from Claim 1. 

If $s_0<m$, define
$$s_1:=\min \Big\{ m,  \min_j \Big(\gamma_j h_j(s_0)\Big)\Big\}.$$
Note that $h_j(s_0)  \gamma_j\ge e^{c_j(m-s_0)}s_0>s_0$ for all $j$, thus $s_1>s_0$. In this setting, we prove:

\med

{\it Claim 2}.  $\displaystyle \min_j\min_{t\in[\tau,2\tau]} x_j(t)\ge s_1.$

\med

Otherwise, there are $t_1\in [\tau, 2\tau]$ and $i\in\{ 1,\dots,n\}$ such that $x_i(t_1)<s_1$ and $x_j(t)\ge x_i(t_1)$ for all $t\in [\tau,t_1]$ and $j\in\{ 1,\dots,n\}$, so (2.11) holds with $t_0$ replaced by $t_1$.
Since $x_i(t_1-\tau_{ik})\ge \min \{s_0,x_i(t_1)\}$, we have 
$h_i(x_i(t_1-\tau_{ik}))\ge \min \{ h_i(x_i(t_1)), h_i(s_0)\}$. We now consider two cases separately.

If $h_i(s_0)\ge h_i(x_i(t_1))$, then $s_0\ge x_i(t_1)$ and we get
$$\eqalign{
0&\ge \left(-d_i+\sum_{j=1}^n a_{ij}\right) x_i(t_1)+\be_i h_i(x_i(t_1)) = \left (-d_i+\sum_{j=1}^n a_{ij}+\be_ie^{-c_ix_i(t_1)}\right ) x_i(t_1) \cr
&> \left (-d_i+\sum_{j=1}^n a_{ij}+\be_ie^{-c_im}\right )x_i(t_1)\ge 0,\cr}
$$
with is not possible.

If $h_i(s_0)< h_i(x_i(t_1))$, then $s_0< x_i(t_1)$. Since $x_i(t_1)<s_1\le \gamma_ih_i(s_0)$, we derive
$$0\ge \big(-d_i+\sum_{j=1}^n a_{ij}\big) x_i(t_1)+\be_i h_i(s_0) > \big (-d_i+\sum_{j=1}^n a_{ij}\big)\gamma_i h_i(s_0)+\be_i h_i(s_0)> 0,$$
which is again a contradiction, ending the proof of Claim 2.

\med

Next, we define by recurrence the sequence
$$s_{k+1}=\min \Big\{ m,  \min_j \Big(\gamma_j h_j(s_k)\Big)\Big\}.$$
If $s_k=m$ for some $k\ge 0$, then 
$\gamma_j h_j(s_k)=\gamma_j e^{-c_jm} m\ge m,$
 hence $s_p=m$ for all $p>k$. In this case, the result follows from Claim 1. Otherwise, 
 $$s_{k+1}= \min_j \Big(\gamma_j h_j(s_k)\Big)\ge \min_j e^{c_j(m-s_k)}s_k>s_k,
 \eq(2.12)$$ and $(s_k)$ is strictly increasing. For $s^*=\lim s_k$,   from (2.12) we have
 $$0<s^*\le m \q {\rm and}\q s^*\ge \min_j e^{c_j(m-s^*)}s^*\ge s^*,$$
 and therefore $s^*=m$. On the other hand,  Claim 2 and an inductive argument imply that 
$\displaystyle \min_j\min_{t\in[k\tau,(k+1)\tau]} x_j(t)\ge s_k, k\ge 0$, and we get $\liminf_{t\to\infty} x_j(t)\ge s^*=m$ for $1\le j\le n$.
\ter

\smal

The permanence of (1.1) is  now an immediate consequence of the lemma above.

\proclaim {Theorem 2.3}. If  \vskip 2mm
{\bf (A1')} $\displaystyle  \exists \ c=(c_1,\dots,c_n)>0:\ {{\be_ic_i}\over {d_ic_i-\sum_{j=1}^n a_{ij}c_j}}>1,\q i=1,\dots,n,$\vskip 2mm
\noindent holds, then  system (1.1) is  uniformly persistent, and thus permanent.

{\it Proof}. The changes of variables $x_i\mapsto c_i^{-1}x_i=\bar x_i, 1\le i\le n,$ transform (1.1) into 
$$\bar x_i'(t)=-d_i\bar x_i(t)+\sum_{j=1}^n \bar a_{ij}\bar x_j(t)+\sum_{k=1}^m \be_{ik}    \bar x_i(t-\tau_{ik})e^{-c_i\bar x_i(t-\tau_{ik})},\q i=1,\dots,n,$$
where $\bar a_{ij}={{c_j}\over {c_i}}  a_{ij}$.
After dropping the bars, we may consider system (2.8), for which condition {\bf (A1)} is satisfied.

 Condition {\bf (A1)} is equivalent to $\be_i>d_i-\sum_{j=1}^n a_{ij}>0$ for $i=1,\dots,n$. 
Choose $L>\max_i (c_i^{-1})$ with $L\ge (\max_i \gamma_i)e^{-1}$ and $m\in (0,c_i^{-1})$ with $ m\le \min_i(c_i^{-1}\log \gamma_i)$.
For $\vare>0$ fixed, let $L_\vare=L+\vare $  and $m_\vare\in (0,m)$ such that $h_i(m_\vare)\le h_i(L_\vare)$. For any positive solution $x(t)$ of  (2.8), let $u_i=\limsup_{t\to\infty}x_i(t)$ and $v_i=\liminf_{t\to\infty}x_i(t)$. Note that $\max_{x\ge 0}h_i(x)=e^{-1}$ for $1\le i\le n$. As in the proof of Theorem 2.1, from (2.1) we deduce that  $\max_i u_i\le \gamma_i e^{-1}<L_\vare$. From Lemma 2.1, we now have $ \min_iv_i>m_\vare$.  By letting $\vare\to 0^+$, we obtain 
$$m\le\liminf_{t\to\infty} x_i(t;\var)\le \limsup_{t\to\infty} x_i(t;\var)\le L,\q 1\le i\le n,
$$ 
for all solutions $x(t;\var)$ of  (2.8) with initial condition $\var\in C^+_0$. For positive solutions of (1.1), we therefore obtain
$$c_im\le\liminf_{t\to\infty} x_i(t;\var)\le \limsup_{t\to\infty} x_i(t;\var)\le c_iL,\q 1\le i\le n. 
\eq(2.13)$$\ter

%

{\bf Remark 2.1}. Consider  a general system (1.1) with coefficients $d_i$  positive, but not given by (1.3). Clearly,
Theorem 2.1 remains true under the additional  condition of $D-A$ being a non-singular M-matrix;  and Theorem 2.3 is valid without further assumptions, since  {\bf (A1')}  implies in particular that  $D-A$ is a non-singular M-matrix, because
$(D-A)c>0$ for some vector $c>0$ [5].

 \med
 
Rather than the estimates (2.13), one can actually give explicit uniform lower and upper bounds for solutions of (1.1),  if lower and upper bounds for the coefficients $\gamma_i$ as defined in (2.9) are known.

\proclaim {Theorem 2.4}.
Assume $e^\al\le \gamma_i\le e^\be, i=1,\dots,n$, with $0<\al<\be, \be >1$.
Then any positive solution $x(t)=(x_1(t),\dots,x_n(t))$ of (1.1) satisfies
$$\min\{\al, \exp \left (\al+\be-1-e^{\be -1}\right)\} \le \liminf_{t\to\infty} x_i(t)\le \limsup_{t\to\infty} x_i(t)\le e^{\be -1},\q i=1,\dots,n.
$$

{\it Proof}. As before, we define $h(x)=xe^{-x}$ for $x\ge 0$. 
If $\max_j u_j=u_i$ for some $i$, from Theorem 2.1 we obtain $(d_i-\sum_{j=1}^n a_{ij}) u_i\le d_iu_i-\sum_{j=1}^n a_{ij} u_j\le \be_i e^{-1}$, which yields $u_i\le \gamma_i e^{-1}\le e^{\be-1}$. Since
$e^{\be-1}>1$, from Lemma 2.1 with $c_1=\cdots=c_n=1$, we have $v_i\ge m, 1\le i\le n$, where $m\in (0,1)$ and is such that $m\le \al$ and $h(m)\le h(e^{\be-1})$.

We now argue as in the proof of Theorem 2.1. Fix a small $\vare>0$, and $T\ge 0$ such that $m-\vare\le v_i-\vare\le x_i(t)\le e^{\be-1}+\vare$ for $t\ge T$ and $1\le i\le n$. Without loss of generality, take $T=0$.
For an arbitrary $t>0$,
$x_i'(t)\ge -d_ix_i(t)+\sum_{j=1}^n a_{ij}(v_j-\vare)+\be_i \min \{h(m-\vare), h(e^{\be-1}+\vare)\}$, implying that
$e^{d_it}x_i(t)\ge x_i(0)+(e^{d_it}-1) \eta_i (\vare)/d_i, t\ge 0$, where  $ \eta_i (\vare)=\sum_{j=1}^n a_{ij}(v_j-\vare)+\be_i \min \{h(m-\vare), h(e^{\be-1}+\vare)\}$.  
Hence we obtain 
$$x_i(t)\ge x_{i}(0)e^{-d_it}+d_i^{-1} \eta_i(\vare) (1-e^{-d_it}),\q i=1,\dots,n.$$
By letting $\vare\to 0^+$ and $t\to\infty$, this leads to
   $v_i\ge d_i^{-1}(\sum_{j=1}^n a_{ij}v_j+\be_ih(m)),$
   for $i=1,\dots,n$.
   For $v_k=\min_i v_i$, this inequality yields 
   $$v_k\ge \gamma_k h(m)\ge e^\al h(m)=e^\al \min\{ h(\al),h(e^{\be -1})\}=\min\{\al, \exp \left (\al+\be-1-e^{\be -1}\right)\}.$$
\vskip -3mm 
 \ter
   
 \med
 
In spite of  the explicit estimates provided by Theorem 2.4, clearly the criterion for the uniform persistence in Theorem 2.3  is more general.
\med

{\bf Example 2.2}.  In (1.1), let $n=2,m=1,\beta_{1}=1,\beta_{2}=3,a_{12}=a_{21}=1,d_{1}=3,d_{2}=2.$
Then  $M=\pmatrix{
-2 & 1\cr
1 & 1\cr}$
 and  $\gamma_{1}<1$, hence  {\bf (A1)} is not satisfied, so Theorem 2.4 does not apply directly. 
 However, it is easy to check that  hypothesis {\bf (A1')}  is satisfied with  any $c_1,c_2>0$ such that $2c_1<c_2<3c_1$,
and therefore we are able to conclude that system (1.1) is permanent.

\bigskip

{
\input epsf
\centerline{\epsfxsize=6cm
(a) \epsfbox{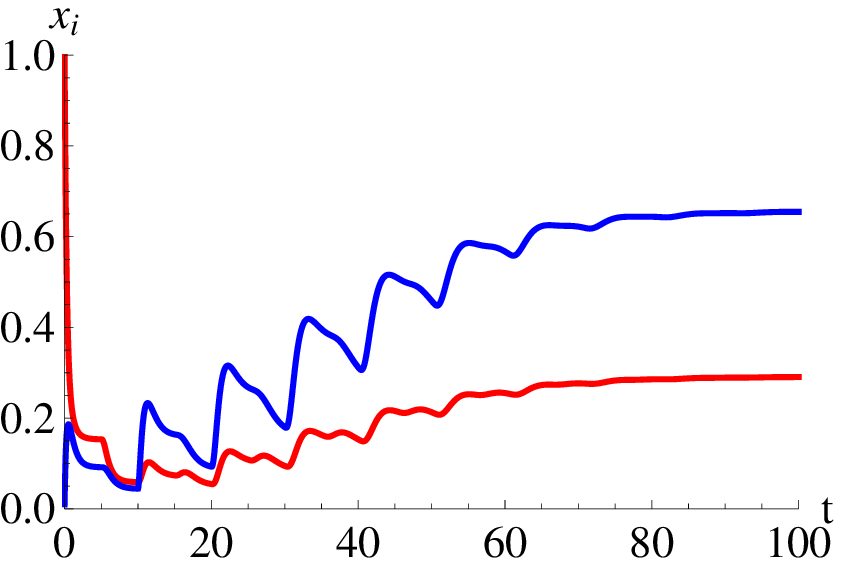}
\epsfxsize=6cm
\hskip0.2in (b) \epsfbox{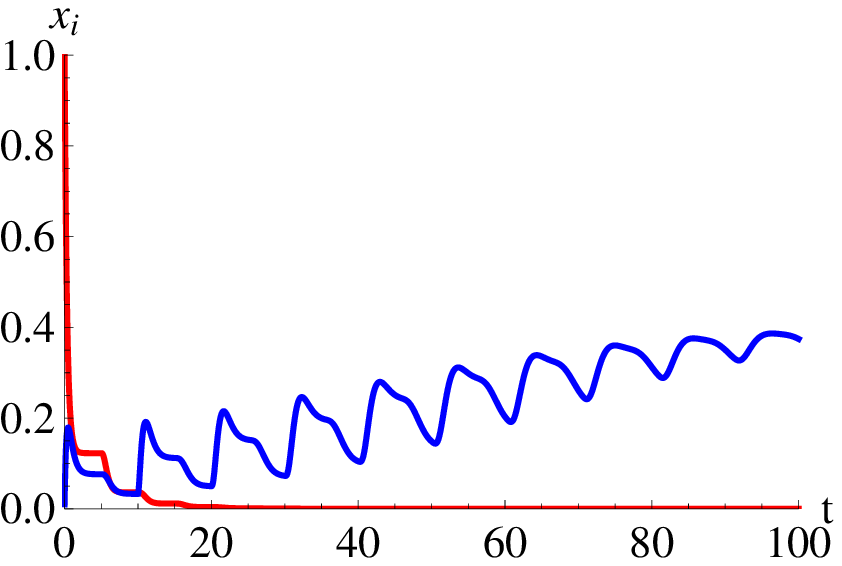}}
\noindent {Figure 1. In (a), Example 2.2 is depicted with $\tau_1=5$ and $\tau_2=10$. (A1) is not satisfied, but (A1') is, and also $s(M)> 0$, hence the population persists on both patches. Furthermore, one can check that the conditions of Theorem 5.2 hold and the positive equilibrium is globally asymptotically stable. In (b), we set $a_{12}=0$, other parameters are the same. Then $s(M)=1>0$, but we are in the reducible case of Example 2.1, and the population becomes extinct on the first patch.}
}

\med

 Next result establishes that  $s(M)>0$ is a  criterion for the uniform persistence of the total population, i.e., the uniform $\rho$-persistence of (1.1) in the sense of Smith and Thieme's nomenclature [16] with $\rho(\phi)=\sum_{i=1}^n \phi_i(0)$; moreover,  in the case of an irreducible matrix $A$, the persistence is uniform in all patches. 
 It will be shown in the next section that this criterion is sharp.
  In the theorem below, we use the norm $|x|=\sum_{i=1}^{n}|x_{i}|$ in $\R^{n}$, so $\rho(\phi)=|\phi(0)|$ for all $\phi \in C^+_0$.


\proclaim{Theorem 2.5}.  Assume $s(M)>0$. Then for system (1.1) the total
population strongly uniformly persists. If $M$ is irreducible, then the population strongly uniformly persist on each patch. If $M$ is reducible,
there exists at least one patch, where the population strongly uniformly persists.

{\it Proof}. The proof is organized in three steps.

\medskip

{\it (i) Finding an irreducible block with positive spectral bound}

\medskip

If $M$ is reducible, then (after a permutation of the variables), it can be written in the diagonal form $$ M=\pmatrix{
M_{11}&\dots  &M_{1\ell}\cr
{}&\ddots&{}\cr
0&\dots& M_{\ell\ell}\cr},
$$
where $ M_{lm}$ are $n_l\times n_m$
matrices, with  $M_{ll}$  irreducible  $n_l\times n_l$ blocks,
$\sum_{l=1}^\ell n_l=n$. Then $s(M)=\max\{ s(M_{\ell \ell}):{i=1,\dots,\ell}\}$,
and there exists an index $\kappa \leq \ell$ such that $s(M_{\kappa\kappa})>0$. Let $\underline \kappa:=\sum_{l=1}^{\kappa-1} n_l+1$ and $\overline \kappa:=\sum_{l=1}^{\kappa} n_l$. Define the index set $\Omega:=\{i \in {\rm I\!N} : \underline \kappa \leq i \leq  \overline \kappa \}$, then $|\Omega|=n_\kappa>0$. Now consider the following subsystem of (1.1), which corresponds to the $\kappa$th block:
$$x_i'(t)=-d_ix_i(t)+\sum_{j\in \Omega} a_{ij}x_j(t)+\sum_{j\notin \Omega} a_{ij}x_j(t)+\sum_{k=1}^m \be_{ik}    x_i(t-\tau_{ik})e^{-x_i(t-\tau_{ik})},\q i \in \Omega. \eq(2.14)$$
In the sequel we let $p_i(t):=\sum_{j\notin \Omega} a_{ij}x_j(t) \geq 0$ for all $i \in \Omega$, and let $\rho^\kappa (\phi):=\sum_{j \in \Omega} \phi_j(0)$. We use the notation $M_{\kappa\kappa}=A_{\kappa\kappa}+B_{\kappa\kappa}-D_{\kappa\kappa}$, where $A_{\kappa\kappa},B_{\kappa\kappa},D_{\kappa\kappa}$ are $n_\kappa\times n_\kappa$ matrices, corresponding to the $\kappa$th block in $A,B,D$. If $M$ is irreducible, we have only one block $M_{11}=M$, and in this case $|\Omega|=n$ and $p_i(t)=0$ for all $i=1,\dots,n$.

\medskip

{\it (ii) Uniform weak persistence of the total population of an irreducible block with positive spectral bound}

\medskip

Consider (2.14). For any $0<\varepsilon<1,$ we define the auxiliary system 
$$
w_{i}^{\prime}=-d_{i}w_{i}(t)+\sum_{k=1}^{m}\beta_{ik}(1-\varepsilon)w_{i}(t-\tau_{ik})+\sum_{j\in \Omega}a_{ij}w_{j}(t),
 \quad i \in \Omega,\eq(2.15)$$
and the auxiliary matrix $M_{\kappa\kappa}({\varepsilon})=A_{\kappa\kappa}+B_{\kappa\kappa}({\varepsilon})-D_{\kappa\kappa},$ where
$$
B_{\kappa\kappa}({\varepsilon})=diag(\beta_{\underline \kappa}(1-\varepsilon),\beta_{\underline \kappa+1}(1-\varepsilon),...,\beta_{\overline \kappa-1}(1-\varepsilon),\beta_{\overline \kappa}(1-\varepsilon)).
$$
If $s(M_{\kappa\kappa})>0,$ then also $s(M_{\kappa\kappa}(\varepsilon))>0$ for sufficiently small
$\varepsilon.$ Fix such an $\varepsilon.$ Since $M_{\kappa\kappa}(\varepsilon)$
(and thus also $M_{\kappa\kappa}(\varepsilon)^{T}$) is a cooperative irreducible
matrix, $s(M_{\kappa\kappa}(\varepsilon))$ is a simple dominant eigenvalue with
a positive eigenvector. Let $q$ be the positive vector that corresponds
to the transpose of $M_{\kappa\kappa}({\varepsilon}),$ i.e. $M_{\kappa\kappa}(\varepsilon)^{T}q=s(M_{\kappa\kappa}(\varepsilon))q.$

Define for any positive solution segment $w_{t}$ of system (2.15)
the vector $y(t)$ by 
$$
y_{i}(t)=\left(w_{i}(t)+\sum_{k=1}^{m}\beta_{ik}\left(1-\varepsilon\right)\int_{t-\tau_{ik}}^{t}w_{i}(u)du\right).
$$
We construct the Lyapunov functional $V:=\langle y(t),q\rangle$ (here
$\langle\cdot,\cdot\rangle$ denotes the Euclidean scalar product). Then it
is easily seen that $y(t)$ satisfies the relation 
$$
y^{\prime}(t)=M_{\kappa\kappa}(\varepsilon)w(t)
$$
and we have 
$$
{{dV(t)}\over{dt}}=\langle y^{\prime}(t),q\rangle=\langle M_{\kappa\kappa}(\varepsilon)w(t),q\rangle=\langle w(t),M_{\kappa\kappa}(\varepsilon)^{T}q\rangle=\langle w(t),s(M_{\kappa\kappa}(\varepsilon))q\rangle>0,
\eq(2.16)$$
because in the last scalar product all terms are positive. Hence
$V$ is increasing and $V>0$ except at zero, so either $\lim_{t\to\infty}V(t)=\infty$
or $\lim_{t\to\infty}V(t)=V_{*}<\infty$ with $V_{*}>0$. We claim
that the latter case is not possible. Suppose the contrary: then by
the fluctuation lemma there is a sequence $t_{l}\to\infty$ as $l\to\infty$
such that $V(t_{l})\to V_{*}$ and $V'(t_{l})\to0$. Then from (2.16)
it follows that $w(t_{l})\to0$. Given that $w_{i}'(t)\geq-d_{i}w_{i}(t)$,
we have that $w_{i}(s)\leq e^{d_{i}\tau}w_{i}(t)$ for any $s\in[t-\tau,t]$,
consequently $y_{i}(t)\leq w_{i}(t)(1+\beta_{i}(1-\varepsilon)\tau e^{d_{i}\tau}).$
As $w(t_{l})\to0$, necessarily $y(t_{l})\to0$ and thus $V(t_{l})\to0$
which is a contradiction. Thus, only $\lim_{t\to\infty}V(t)=\infty$
is possible.

Now consider a positive solution $x(t)$ of (1.1), and let $\tilde x(t)=(x_{\underline \kappa}(t),\dots,x_{\overline \kappa}(t))^T$. There
is a $\delta_{0}=\delta_{0}(\varepsilon)>0$ such that $e^{-\xi}>\left(1-\varepsilon\right)$
for $\xi\in[0,\delta_{0}].$ Then $\beta_{ik}\xi e^{-\xi}\geq\beta_{ik}\xi\left(1-\varepsilon\right)$
for all $i=1,...n,\, k=1,\dots,m$ and $\xi\in[0,\delta_{0}].$ Define
the set $U_{\varepsilon}$ by 
$$
U_{\varepsilon}=\{\psi\in C^{+}([\tau,0],\R^{|\Omega|}):||\psi_{i}||\le \delta_{0}\ {\rm for\ all \ }i \in \Omega\}.
$$
Suppose that there is a $t_{0}$ such that $\tilde x_{t}\in U_{\varepsilon}$
for all $t\geq t_{0}$. Then we can consider a solution $w(t)$ of
(2.15) for $t\geq t_{0}$ with $w_{t_{0}}=\tilde x_{t_{0}}$, and
by a standard comparison principle (using $p_i(t)\geq 0$ and $\beta_{ik}\xi e^{-\xi}\geq\beta_{ik}\xi\left(1-\varepsilon\right)$) we obtain $\tilde x(t)\geq w(t)$ for all
$t\geq t_{0}$, and $\tilde x_{t}\in U_{\varepsilon}$ implies $w_{t}\in U_{\varepsilon}$
for all $t\geq t_{0}$, which contradicts $V(t)\to\infty$.

Therefore, there is a sequence $t_{l}\to\infty$ as $l\to\infty$
such that $\tilde x_{t_{l}}\notin U_{\varepsilon}$. Then for each $t_{l}$
there is a $j(l) \in \Omega$ such that $||(\tilde x_{t_{l}})_{j(l)}||>\delta_{0}$,
thus there is a $t_{l}^{*}\in[t_{l}-\tau,t_{l}]$ such that $\tilde x_{j(l)}(t_{l}^{*})>\delta_{0}$.
By $\tilde x'_{j(l)}(t)\geq-d_{j(l)}\tilde x_{j(l)}(t)$ we have $\tilde x_{j(l)}(t_{l})\ge \tilde x_{j(l)}(t_{l}^{*})e^{-d_{j(l)}(t_{l}-t_{l)}^{*})}\ge e^{-d_{j(l)}\tau}\delta_{0}$
, thus 
$$
|\tilde x(t_{l})|\geq\delta:=\min_{i=1,\dots,n}\{e^{-\tau d_{i}}\delta_{0}\},
$$
and we obtain that $\limsup_{t\to\infty}|\tilde x(t)|\geq\delta,$ hence
we obtain the uniform weak persistence of the total population on the patches of the $\kappa$th block.

We conclude that system (1.1) is uniformly weakly $\rho^\kappa$-persistent
with $\rho^\kappa(\phi)=\sum_{i\in \Omega}\phi_{i}(0)$, which represents the
persistence of the total population of the patches of the $\kappa$th block.
\med

{\it (iii) Uniform strong persistence on each patch of an irreducible block with positive spectral bound}

\med

To show the uniform strong $\rho^\kappa$-persistence (i.e. there is a $\theta>0$
such that $\liminf_{t\rightarrow\infty}\rho^\kappa(x_{t})>\theta),$ we can
apply Theorem 4.5 of [16, Chapter 4.1].
By the dissipativity (Theorem 2.1), there exists a compact global attractor of system (1.1)
(by [7], Theorem 3.4.8), and the conditions of Theorem 4.5
of [16] hold, which guarantees the uniform strong
$\rho^\kappa$-persistence. Next we show the persistence of the population
in each patch of the $\kappa$th block. We shall use the persistence functions $\rho_{i}(x_{t})=x_{i}(t),$
which express the actual population on patch $i$. Let $\epsilon\in(0,\theta)$,
where $\theta$ corresponds to $\rho^\kappa$-persistence, i.e. $\liminf_{t\rightarrow\infty}\rho^\kappa(x_{t})>\theta$.
Then for any solution $x_{t}$ there is a sequence $t_{l}\to\infty$
as $l\to\infty$ such that $\sum_{i\in \Omega}x_{i}(t_{l})>\theta-\epsilon$
for all $l$. Then there must be an index $j \in \Omega$ such that $x_{j}(t_{l})>{{\theta-\epsilon}\over n}$
holds for infinitely many $t_{l}$. We may assume $j=\underline \kappa$. Thus, $\limsup_{t\to\infty}x_{\underline \kappa}(t)\geq{{\theta-\epsilon}\over n}$,
and the system is uniformly weakly $\rho_{\underline \kappa}$-persistent. We can
apply again Theorem 4.5 of [16] to conclude the uniform
strong $\rho_{\underline \kappa}$-persistence, thus there is an $\eta_{\underline \kappa}>0$ such
that $\liminf_{t\to\infty}x_{\underline \kappa}(t)>\eta_{\underline \kappa}$ and the population persists
on patch $\underline \kappa$. By the irreducibility of $M_{\kappa\kappa}$, there is an index $j \in \Omega$, such
that $a_{j\underline \kappa}>0$. We may assume $j=\underline \kappa+1$, then $x'_{\underline \kappa+1}(t)\geq-d_{\underline \kappa+1}x_{\underline \kappa+1}(t)+a_{\underline \kappa+1,1}x_{\underline \kappa}(t)$,
thus $\liminf_{t\to\infty}x_{\underline \kappa+1}(t)>\eta_{\underline \kappa+1}$, where we can choose
$\eta_{\underline \kappa+1}=\eta_{\underline \kappa}a_{\underline \kappa+1,\underline \kappa}/d_{\underline \kappa+1}$. By the irreducibility of this block, we can reach
all patches inductively and by choosing $\eta=\min_{i\in \Omega}\{\eta_{i}\}$
we have proved the statement of the theorem, and the population strongly
uniformly persists on each single patch $i \in \Omega$. In the irreducible case, 
$\Omega$ contains all indices $i=1,\dots,n$ and the population strongly
uniformly persists on each patch.
\ter

\bigskip

{
\input epsf
\centerline{\epsfxsize=6cm
(a) \epsfbox{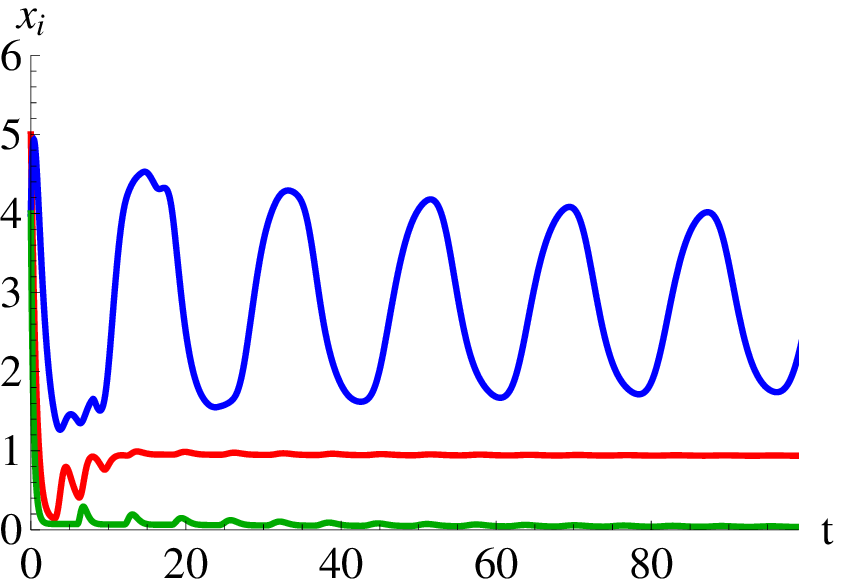}
\epsfxsize=6cm
\hskip0.2in (b) \epsfbox{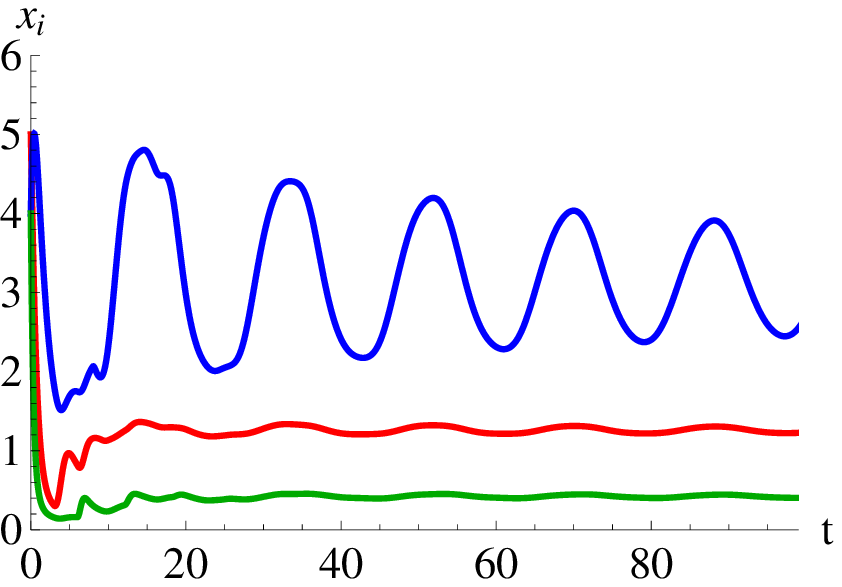}}
\noindent {Figure 2. Illustration of a system with three patches. In (a), parameters are set to $n=3$, $m=1$, $\beta_1=5$,$\beta_2=10$, $\beta_3=3$, $d_1=2$, $d_2=1$, $d_3=3$, $a_{12}=a_{31}=a_{32}=0$, $a_{13}=a_{21}=a_{23}=1$, $\tau_1=3$, $\tau_2=8$, $\tau_3=6$. Then $M$ is reducible but $s(M)=9>0$. We can observe different behavior on the patches: oscillation, convergence to a positive value, extinction. In (b), parameters are the same, except that $a_{12}=a_{31}=a_{32}=0.1$, thus $M$ is irreducible and the system is persistent.}
}

\

{\bf 3. Extinction}

\med

In this section, a sharp criterion for  the global asymptotic stability of the trivial equilibrium of (1.1) is established. In biological terms, this means the extinction of the population in all patches.

\proclaim{Theorem 3.1}.  Suppose that  $s(M)\le 0$. Then the equilibrium 0 of (1.1) is GAS (in $C^+$).

{\it Proof}. If $s(M) <0$, or if $s(M)=0$ and $A=[a_{ij}]$ is an irreducible matrix, the global asymptotic stability of $x=0$ follows from Theorems 2.1 and 3.1 in [3], respectively; for the latter case, the framework in [19] was used.

 Now, suppose that $A$ is reducible and $s(M)=0$. After a permutation of the variables in (1.1),  we may suppose that $A$ has the form
$$ A=\pmatrix{
A_{11}&\dots  &A_{1\ell}\cr
{}&\ddots&{}\cr
0&\dots& A_{\ell\ell}\cr},
$$
where $ A_{km}$ are $n_k\times n_m$
matrices, with  $A_{kk}$  irreducible  $n_k\times n_k$ blocks,
$\sum_{k=1}^\ell n_k=n$.  (According to our definition, here a square matrix of size one is always irreducible; cf.~e.g.~Appendix A of [17].)

We prove the result for $\ell=2$; the general case follows by induction.  Suppose that $n_1+n_2=n$ and $a_{ij}=0$ for $n_1+1\le i\le n ,1\le j\le n_1$, so that
$$A=\pmatrix {A_{11}&A_{12}\cr 0&A_{22}\cr},\q  M=\pmatrix {M_{11}& M_{12}\cr 0& M_{22}\cr},
\eq(3.1)$$
where $A_{ij}, M_{ij}$ are $n_i\times n_j$ blocks and $M_{ii}$ are  irreducible matrices.  Since $\sigma (M)=\sigma (M_{11})\cup \sigma (M_{22})$, we have $s (M_{ii})\le 0, i=1,2$.  

Write a solution $x(t)=x(t;\var)$ (for $\var\in C^+$) of (1.1) as $x(t)=(y(t),z(t))\in \R^{n_1}\times \R^{n_2}$ according to the decomposition of $M$ in (3.1). The result for the irreducible case implies that 0 is the unique equilibrium of (1.1), and that $z(t)\to 0$ as $t\to\infty$. If suffices to show that $y(t)\to 0$ as $t\to\infty$.

Since $s(M_{11})\le 0$, then  $-M_{11}$ is an  M-matrix; moreover, $-M_{11}$ is an irreducible matrix, therefore that  there exists a positive vector $c=(c_1,\dots,c_{n_1})$ such that $M_{11}c\le 0$ [5], i.e.,
$$\be_i-d_i+\sum_{j=1}^{n_1} {c_j\over c_i}a_{ij}\le 0,\q i=1,\dots,n_1.
\eq(3.2)$$
Rewrite system (1.1) with the change of variables $\bar y_i=c_i^{-1}y_i, i=1,\dots,n_1$.
Dropping the bars for the sake of simplification, we get
$$\eqalign{
y_i'(t)&=-d_iy_i(t)+\sum_{j=1}^{n_1} {c_j\over c_i}a_{ij}y_j(t)
+\sum_{k=1}^m \be_{ik}    y_i(t-\tau_{ik})e^{-c_iy_i(t-\tau_{ik})}+g_i(t),\q i=1,\dots,n_1\cr
z_p'(t)&= -d_pz_p(t)+\sum_{j=n_1+1}^{n} a_{pj}z_p(t)
+\sum_{k=1}^m \be_{pk}    z_p(t-\tau_{pk})e^{-z_p(t-\tau_{pk})}, \q p=1,\dots,n_2
\cr},\eq(3.3)$$ 
where $g_i(t):=\sum_{k=1}^{n_2} a_{i(n_1+k)}z_k(t)\to 0$ as $t\to\infty$. 
Next, define 
$$u_j=\limsup_{t\to\infty} y_j(t),\eq(3.4)$$
 where $y_j,z_p$ satisfy (3.3).
We need to prove that $u:=\max_{1\le j\le n_1} u_j=0$.

Suppose  that $u>0$. For each $i\in\{1,\dots,n_1\}$ such that $u_i=u$, by the fluctuation lemma  there is a  sequence $(t_k)$, with $t_k\to\infty$, 
$y_i(t_k)\to u_i, y_i'(t_k)\to 0$. Choose $\vare \in (0,u_i)$. For $t$ and $k$ large,  we have  $y_i(t_k)\ge u_i-\vare$, $y_j(t)\le u_j+\vare,  j=1,\dots,n_1$, and $0\le g_i(t)\le \vare $, leading to
$$y_i'(t_k)\le-d_i(u_i-\vare)+\sum_{j=1}^{n_1} {c_j\over c_i}a_{ij}(u_j+\vare)
+ \be_i   (u_i+\vare)+\vare.$$ 
By letting $\vare\to 0^+$ and $k\to\infty$,  from (3.2) we get
$$0\le (\be_i-d_i)u_i+\sum_{j=1}^{n_1} {c_j\over c_i}a_{ij}u_j
\le \Big(\be_i -d_i+\sum_{j=1}^{n_1} {c_j\over c_i}a_{ij}\Big )u_i \le  0.
\eq(3.5)$$
This leads to
$$\be_i-d_i+\sum_{j=1}^{n_1} {c_j\over c_i}a_{ij}= 0,\q \sum_{j=1}^{n_1} {c_j\over c_i}a_{ij}(u_j-u_i)=0, \q {\rm if}\q u_i=u.
\eq (3.6)$$
 On the other hand, reasoning as in the proof of Theorem 2.1, and since $\lim_{t\to\infty}z_p(t)=0$ for $1\le p\le n_2$,  (3.2) and (3.6) yield   the estimate
$$\be_iu=d_iu-\sum_{j=1}^{n_1} {c_j\over c_i}a_{ij}u\le \be_i(c_ie)^{-1},$$
implying that $u\le (c_ie)^{-1}$. In particular, for any $\vare >0$ and $i$ such that $u_i=u$,  the bounds
$0\le y_i(t)<(u+\vare)<1/c_i$ hold
 for $t>0$ large.

Next, for  $i$ such that $u_i=u$ consider again a  sequence  $(t_k)$ as above.
 For $\vare >0$ small and $k$ large, 
$$
 y_i'(t_k)\le-d_i(u-\vare)+\sum_{j=1}^{n_1} {c_j\over c_i}a_{ij}(u+\vare)
 +\sum_{q=1}^m \be_{iq} h_i(y_i(t_k-\tau_{iq}))+\vare,
 $$ 
 where  $h_i(x)=xe^{-c_ix}$. The  functions $h_i$ are strictly increasing for $0\le x\le 1/c_i$, hence $h_i(y_i(t_k-\tau_{iq}))\le h_i(u+\vare)$ for $k$ large.  From (3.6), and  letting $\vare\to 0^+$ and $k\to\infty$, we thus obtain
$$0\le  \be_iu\, (e^{-c_i u}-1)<0,$$
which is not possible. This shows that $u=0$, and  the proof is complete.\ter 

\med

In view of Theorems 2.2, 2.5 and 3.1, we therefore have a sharp threshold criterion  for {\it extinction} versus {\it uniform persistence} of the {\it total population} in the general case; and in the case of an irreducible matrix $A$, we have a sharp threshold criterion for extinction versus uniform persistence of the population in {\it all patches}. Such consequences are formulated in the following two theorems.

\proclaim{Theorem 3.2}.  If  $s(M)\le 0$, the equilibrium 0 of (1.1) is GAS; while if $s(M)>0$, the total population is uniformly persistent. 

\proclaim{Theorem 3.3}.  Suppose that  the matrix   $A$ is irreducible.  Then:
(i)   if  $s(M)\le 0$, the equilibrium 0 of (1.1) is GAS; (ii) if $s(M)>0$, system (1.1) is uniformly persistent, i.e., the population uniformly persists on each patch.  Moreover, $s(M)>0$ if and only if  there is a positive vector $c\in\R^n$ such that $Mc>0$.

As observed,   $s(M)>0$ is a sharp condition for the uniform persistence of  (1.1) in the irreducible case, whereas this criterion fails  in the case of reducible  community matrices. In the latter case, while the total population uniformly persists if $s(M)>0$, the population can become extinct on some of the patches (see Example 2.1). However, by Theorem 2.3  the uniform persistence follows under the stronger hypothesis {\bf (A1')}.

Two  final notes in this section open the present framework to possible generalizations.


\med

{\bf Remark 3.1.} Theorem 3.1 is also valid for a system (1.1) without condition (1.3). In fact, since $s(M)  \le 0$ is equivalent  to saying that $-M=D-A-B$ is an M-matrix, and 
  $\underline\be =\min_i \be_i$ is strictly positive, then $s(M)  \le 0$ implies that $D-A\ge M+\underline\be I$ is a non-singular M-matrix [5, Theorem 5.3]. In view of this, by Theorem 2.1 and Remark 2.1,  all solutions of (1.1) with initial conditions in $C^+$ are bounded, and in this way the limits in (3.4) are well-defined.
  
  \med
  
{\bf Remark 3.2.}  Some results in Sections 2 and  3  can be  extended in a natural way to a more general class of delayed systems with patch structure of the form $x'_i(t)=-d_{i}x_{i}(t)+\sum_{j=1}^{n}a_{ij}x_j(t)+b_i(x_{t,i}), \, 1\le i\le n$, where the birth functions $b_i:C([-\tau,0];\R)\to \R_+$ are $C^1$-smooth, bounded, with $b_i(0)=0, Db_i(0)(1)=\be_i$, and satisfy some additional conditions. 
Nevertheless, we emphasize that the uniform estimates provided by Theorems 2.3 and 2.4 are valid  for the specific  Ricker-type non-linearity  only. Also, the main result on the global asymptotic stability of the positive equilibrium, which will be presented in Section 5,   depends heavily on the shape of the non-linearity $h(x)=xe^{-x}$, and cannot be extrapolated for a more general class of population models.

\

{\bf 4. Existence of a positive equilibrium}

\med

Together with (1.1), we consider the ODE model  in the positive cone $\R^n_+$:
$$x_i'=-d_ix_i+\sum_{j=1}^n a_{ij}x_j+\be_i  x_ie^{-x_i}=:f_i(x),\q i=1,\dots,n.\eq(4.1)$$
For all $i\in \{ 1,\dots, n\}$ and $x\in \R^n_+$ with $x_i=0$, we have $f_i(x)\ge 0$, thus the positive cone $\R^n_+$ is positively invariant for (4.1). 

The ODE (4.1) may be seen as the particular case of (1.1) with $\tau=0$. 
Clearly, systems (1.1) and (4.1) share the same  equilibria. In this section, we look for equilibria of (4.1). 

In the following, we adopt  some definitions and notation of Hofbauer [8], namely the definition of a saturated equilibrium (or saturated  fixed point).
For an  ODE system $x'=f(x)$ for which $\R^n_+$ is  positively invariant, if an equilibrium point $x^*$ lies on the frontier of $\R^n_+$, say $x^*=(0,\dots, 0,x_{p+1}^*,\dots, x_n^*)$, then necessarily the Jacobian matrix $Df(x^*)$ has the form (cf. [8])
$$Df(x^*)=\left [\matrix{ C&0\cr D&E \cr}\right ],$$  
where $C$ is a $p\times p$ matrix,  called the {\it external part}  of $Df(x^*)$. 

\med

{\bf Definition 4.1}.  For an ODE system $x'=f(x)$, positively invariant in $\R^n_+$,
an equilibrium $x^*\ge 0$ is said to be a {\bf saturated equilibrium} if $x^*$ is an equilibrium and:
 (i) either $x^*\in int ( \R^n_+)$ and $Df(x^*)$ is stable, i.e., $s\big ( Df(x^*)\big )\le 0$; 
(ii) or  $x^*\in fr( \R^n_+)$, $x^*=(0,\dots, 0,x_{p+1}^*,\dots, x_n^*)$, and
$Df(x^*)=\left [\matrix{ C&0\cr D&E \cr}\right ]$,  where the $p\times p$ matrix $C$  is  stable, i.e.,   $s(C)\le 0$.

An equilibrium $x^*\ge 0$ of (4.1) is said to be {\bf regular} if $\det Df(x^*)\ne 0$; in this case, the {\bf index} of $x^*$ is defined as the sign of $\det (-Df(x^*)).$

\med

With these definitions, note that an asymptotically  stable equilibrium has index $+1$, in any dimension $n$. 

The following theorem  plays an important role in this section.
 
\proclaim{Theorem 4.1}. [8]  Any system
  $x'=f(x)$ for $x\in \R^n_+$, where $f$ is a $C^1$ vector field, which is dissipative and forward  invariant  on $\R^n_+$ has at least one saturated equilibrium; moreover, if all saturated equilibria are regular, the sum of their indices equals $+1$.

For system (4.1), the ODE version of Theorem 2.1 shows that (4.1) dissipative. Consequently, from Hofbauer's theorem we deduce that there is at least a saturated fixed point of (4.1) in the cone $\R^n_+$.

Next, we give sufficient conditions for the existence and stability of a positive equilibrium of (4.1), both for the irreducible and reducible case. A sharp criterion  is obtained when $A$ is  irreducible.

\proclaim{Theorem 4.2}. Assume $A$ is irreducible. If
$s(M)>0$,  there is a unique positive equilibrium $x^*$ of (4.1),  which is GAS in $\R^n_+\setminus \{ 0\}$; if $s(M)\le 0$, zero is a global attractor in $\R^n_+$.

{\it Proof}. The last assertion follows from Theorem 3.1. Now, suppose that $s(M)>0$. From Theorem 4.1, there is a saturated equilibrium of (4.1). Since $A$ is irreducible, 
the Jacobian matrix at an equilibrium $u^*$,
$Df(u^*)=A-D+diag\,  \Big ( \be_ie^{-u_i^*}(1-u_i^*)\Big )_{i=1}^n,$
is also  {irreducible}, thus the only possible saturated equilibrium on the boundary of $\R^n_+$ is zero, for which  the external part of $ Df(0)$ coincides with the full matrix.  
However, condition $s(M)>0$ implies that the linearized equation at 0, $\dot x=Mx$, has an eigenvalue with positive real part, hence zero is an unstable fixed point of (4.1). Consequently,
 there is a positive saturated equilibrium $x^*$. But any other possible positive equilibrium of (4.1) is saturated. In fact, if $u^*>0$ is an equilibrium of (4.1), we have 
$$-Df(u^*)u^*=col\,  \Big ( \be_ie^{-u_i^*}(u_i^*)^2\Big )_{i=1}^n>0.$$
This implies that $-Df(u^*)$ is a non-singular M-matrix (see [5]), which is equivalent to saying that $s(Df(u^*))<0$.
Therefore $u^*$ is regular with index +1. 
Again by  Theorem 4.1 we conclude that the positive equilibrium $x^*$ of (4.1) is unique, and locally asymptotically stable. Since (4.1) is an irreducible and cooperative system, by Theorem 6 of  [8] (see also proof of Lemma 4.2 below)  $x^*$ is a global attractor of all positive solutions $x(t)$. On the other hand,  any solution $x(t)=x(t;x_0)$ of (4.1) with initial condition in $x_0\in \R^n_+\setminus \{ 0\}$ is strictly positive for $t>0$ (cf.~e.g.~
[15]).\ter

\proclaim{Theorem 4.3}. Assume (2.6)  for some $c=(c_1,\dots,c_n)>0$. Then, there is a unique positive equilibrium $x^*$ of (4.1), which is GAS in $int(\R^n_+)$.

{\it Proof}. 
If $A$ is irreducible, (2.6) is equivalent to $s(M)>0$ (cf.~Theorem 3.3).
If $A$ is a reducible matrix, the existence of a globally asymptotically stable positive equilibrium of (4.1) is an immediate consequence of the next two lemmas.\ter

\proclaim{Lemma 4.1}. If (2.6) holds, then there is a unique positive equilibrium of (4.1).

{\it Proof}. As before, write the ODE (4.1) as $x'=f(x)$, for  $f=(f_1,\dots, f_n)$ and $f_i(x)=(\be_ie^{-x_i}-d_i)x_i+\sum a_{ij}x_j$, and designate by $x(t,x_0)$ the solution of (4.1) with initial condition $x(0)=x_0\in\R^n_+$.   For a vector $c$ as in (2.6) , we have $f_i(\vare c)=\vare [-(c_id_i-\sum c_ja_{ij})+c_i\be_i e^{-\vare c_i}]$, hence $f_i(\vare c)>0$ for $\vare>0$ small and $1\le i\le n$. Since (4.1) is cooperative and dissipative, from Corollary 5.2.2 of [15, p.~82], $x(t,\vare c)\to x^*$ as $t\to\infty$ for some $x^*>0$. Clearly $x^*$ is an equilibrium of (4.1).
It suffices to show that  $x^*$  is the unique positive fixed point. 

The case of $A$ irreducible has already been addressed.  Now, suppose that $A$ is reducible, with
$$A=\pmatrix {A_{11}&A_{12}\cr 0&A_{22}\cr},$$
where the $n_i\times n_i$ matrices $A_{ii}$ are  irreducible ,  $i=1,2,\ n_1+n_2=n$. (Recall that this includes the case of  some of the $A_{ii}$ equal to zero  if $n_i=1$.)
The general case where $A$ can be written in a triangular form with $\ell$ irreducible diagonal blocks $A_{ii}$ follows by induction. We write accordingly 
$$M=\pmatrix {M_{11}&M_{12}\cr 0&M_{22}\cr},\q c=\pmatrix{c^{(1)}\cr c^{(2)}\cr},$$
with $n_i\times n_i$ matrices $M_{ii}$ and $c^{(i)}\in\R^{n_i},\ i=1,2$.
Since $Mc>0$, then $M_{22}c^{(2)}>0$, and  Theorem 3.2 yields $s(M_{22})>0$.

For $x(t)=(y(t),z(t))\in \R^{n_1}\times \R^{n_2}$, system (4.1) becomes
$$\eqalignno{
y_i'&=(\be_i   e^{-y_i}-d_i)y_i+\sum_{j=1}^{n_1} a_{ij}y_j +\sum_{k=1}^{n_2} a_{i(n_1+k)}z_k,\q i=1,\dots,n_1&(4.2_a)\cr
z_p'&= (\be_p   e^{-z_p}-d_p)z_p+\sum_{k=1}^{n_2} a_{p(n_1+k)}z_k, \q p=1,\dots,n_2.&(4.2_b)\cr
}$$
Write $x^*=(y^*,z^*)\in \R^{n_1}\times \R^{n_2}$. From the irreducible case, 
 $z^*$ is the unique positive equilibrium of $(4.2_b)$, which is GAS. If $A_{12}=0$, then clearly $y^*$ is the unique positive equilibrium of $(4.2_a)$.
Otherwise, define $l:=A_{12}z^*$ and note that $ l=(l_1,\dots,l_{n_1})\ge 0, l\ne 0$.  Consider the system
$$y_i'=(\be_i   e^{-y_i}-d_i)y_i+\sum_{j=1}^{n_1} a_{ij}y_j +l_i=:g_i(y),\q i=1,\dots,n_1.
\eq(4.3)$$
Obviously 0 is not a fixed point of (4.3). The positive cone $\R^{n_1}_+$   is positively invariant for (4.3). For $u^*=(u_1^*,\dots,u^*_{n_1})$ an equilibrium of (4.3),
$Dg(u^*)=diag\, \Big (\be_ih'(u_i^*)-d_i\Big )+A_{11}$ is irreducible, thus there are no saturated equilibria on the boundary of $\R^{n_1}_+$.
Also, 
$Dg(u^*)u^*=-col\, \Big( (u_i^*)^2e^{-u_i^*}+l_i\Big )_{i=1}^{n_1}<0$, and therefore we conclude that $-Dg(u^*)$ is a non-singular M-matrix, which implies that $u^*$ is regular with index +1.
From Theorem 4.1, we deduce that (4.3) has a unique saturated equilibrium, which is $y^*$. This ends the proof.\ter

\proclaim{Lemma 4.2}. If there exists a unique positive equilibrium $x^*$ of (4.1), then $x^*$ is GAS in $int(\R^n_+)$.

{\it Proof}. Let $x_0\in int(\R^n_+)$. Choose $l,L$, $0<l<1<L$,  such that $lx^*\le x_0\le Lx^*$.  With the same notations as above, we have that
$f_i(lx^*)>lf_i(x^*)=0$ and $f_i(Lx^*)<Lf_i(x^*)=0$. This implies that the components $x_i(t, lx^*)$ are non-decreasing and $x_i(t,Lx^*)$ are non-increasing, for $t\ge 0$ [15, Corollary 5.2.2]. Reasoning as above, let $K_1, K_2$ be such that $x(t,lx^*)\to K_1$ and $x(t,Lx^*)\to K_2$ as $t\to\infty$. Clearly $K_1,K_2$ are positive equilibria, hence $K_1=K_2=x^*$. Since (4.1) is cooperative, $x(t,lx^*)\le x(t,x_0)\le x(t,Lx^*)$, hence $x(t,x_0)\to x^*$ as $t\to\infty$. \ter

\med

The results in Sections 2 to 4 yield some interesting algebraic consequences, which may be useful  in applications.

\proclaim{Theorem 4.4}.  (i) For a cooperative matrix $M$,  if $Mc>0$ for some positive vector $c$, then $s(M)>0$; the converse is true if $M$ is irreducible.\vskip 0cm
(ii) If $M=B-D+A$ for $A,B,D$ as in (1.4), with either (1.3)  or $D-A$  a non-singular M-matrix, then  {\bf (A1')} holds if and only if $Mc>0$ for some positive vector $c$.

{\it Proof}. (i) From Theorems 2.2 and  3.1, condition (2.6) implies $s(M)>0$.  (ii) Obviously, {\bf (A1')} implies (2.6).
If $Mc>0$ for some positive vector $c$, from Theorem 4.3 there is a unique positive equilibrium $x^*>0$ of (4.1) (and (1.1)) (note that the dissipativity of (4.1) follows from $D-A$  being a non-singular M-matrix, in case (1.3) is not satisfied). Consequently, $Bx^*>diag (\beta_i x_i^* e^{-x_i^*})=(D-A)x^*>0$, thus
{\bf (A1')} is satisfied with $c=x^*$.\ter

\

{\bf 5. Global asymptotic stability of the positive equilibrium}

\med

In this section, we give a criterion for the  (absolute) global attractivity of the positive equilibrium.
We shall use an auxiliary result established in  [3].

\proclaim{Lemma 5.1}.  [3] The function $h(x)=xe^{-x}$ satisfies
$$|h(y)-h(x)|< e^{-x} |y-x|\q {\rm for \ all}\ x\in (0,2]\ {\rm and}\ y> 0, y\ne x.
$$

We now prove the main result of this section.

\proclaim {Theorem 5.1}. 
Assume\vskip 0.2cm
{\bf (A2) } $1<\gamma_i\le e^2$,\q  $i=1,\dots,n$, where $\displaystyle \gamma_i:={{\be_i}\over{d_i-\sum_{j=1}^n a_{ij}}}.$
\vskip 0.2cm
\noindent 
Then  the positive equilibrium $ x^*$ for (1.1) is
GAS (in $C_0^+$).

{\it Proof}. Theorems 2.1 and 4.3 guarantee  that all positive solutions of (1.1) are bounded and that there is a unique positive equilibrium $x^*=(x_1^*,\dots,x_n^*)$ of (1.1). For  $x^*_i=\max_j x_j^*$, we obtain $$e^{x_i^*}\le \gamma_i\le e^2,$$ hence $0<x_j^*\le x_i^*\le 2,\, 1\le j\le n$. Thus,
 $x^*$ is locally asymptotically stable (cf. Theorem 2.2 and  [3, Remark 2.1]).

As before,  let $h(x)=xe^{-x}$ for $x\ge 0$, and effect the changes 
$$z_i(t)={{x_i(t)}\over x_i^*}-1,\q 1\le i\le n.\eq(5.1)$$
 System (1.1) becomes
$$z_i'(t)={1\over x_i^*}\left [ -d_ix_i^*z_i(t)+\sum_{j=1}^n a_{ij}x_j^*z_j(t)+\sum_{k=1}^m \be_{ik}    \Big (h(x_i^*+x_i^*z_i(t-\tau_{ik}))-h(x_i^*)\Big )\right],\ i=1,\dots,n.\eq(5.2)$$

 \med
 
 Consider any solution $z(t)=z(t;\phi)$ of (5.2) with initial condition $\phi\in S$, where $S:=\{ \phi=(\phi_1,\dots,\phi_n)\in C([-\tau,0];\R^n): \phi_i(\th)\ge -1$ for $-\tau\le \th <0$ and $\phi_i(0)> -1,\ i=1,\dots,n\}$. 
 Then, there are constants $m,M$, $0<m<M$, with $m-1<z_i(t)<M$ for all $i$ and $t>0$ sufficiently large.
 To prove that $z(t)\to 0$ as $t\to\infty$, we now follow closely some arguments in [3]. 
 
 Fix the maximum norm in $\R^n$, $|x|=\max_{1\le i\le n}|x_i|$ for $x=(x_1,\dots,x_n)$. If $\phi =0$, then $z(t)\equiv 0$. For $\phi\ne 0$, we claim
 that 
 $$|z(t)|< \|\phi\|\ {\rm for}\ t\ge \tau.
 \eq(5.3)$$

 For the sake of contradiction, suppose that (5.3) fails. Then, there exists $T\ge \tau$ such that $|z(T)|\ge \|\phi\|>0$ and $|z(T)|\ge |z(t)|$ for $-\tau \le t\le T$.

Let $i\in\{ 1,\dots,n\}$ be such that $|z(T)|=|z_i(T)|$, and consider the case $z_i(T)>0$ (the case $z_i(T)<0$ is  similar). From the definition of $T$, we have $z_i'(T)\ge 0$. On the other hand, we obtain
$$z_i'(T)= {1\over x_i^*}\left [ -d_ix_i^*z_i(T)+\sum_{j=1}^n a_{ij}x_j^*z_j(T)+\sum_{k=1}^m \be_{ik}    \Big (h(x_i^*+x_i^*z_i(T-\tau_{ik}))-h(x_i^*)\Big )\right].\eq(5.4)$$
Note $T-\tau_{ik}\ge 0$, hence $x_i^*+x_i^*z_i(T-\tau_{ik})$ is strictly positive. By Lemma 5.1, if $z_i(T-\tau_{ik})\ne 0$, then
$$|h(x_i^*+x_i^*z_i(T-\tau_{ik}))-h(x_i^*)|<e^{-x_i^*}x_i^*|z_i(T-\tau_{ik})|\le e^{-x_i^*}x_i^*z_i(T);$$
  and $h(x_i^*+x_i^*z_i(T-\tau_{ik}))-h(x_i^*)=0$ if $z_i(T-\tau_{ik})=0$. Since $\be_i=\sum_k \be_{ik}> 0$, then $\be_{ik}>0$ for some $k$, and   clearly we obtain $\sum_{k=1}^m \be_{ik}    \Big (h(x_i^*+x_i^*z_i(T-\tau_{ik}))-h(x_i^*)\Big )< \be_i e^{-x_i^*}x_i^*z_i(T)$.
 Also,    $|z_j(T)|\le z_i(T)$ for all $j$, and consequently (5.4) yields
$$z_i'(T)<  {1\over x_i^*} \left [ (-d_ix_i^*+\sum_{j=1}^n a_{ij}x_j^*)+\be_i  e^{-x_i^*}  x_i^*\right]z_i(T)=
0,$$
which contradicts the fact $z_i'(T)\ge 0$. 
This proves (5.3). 


\med

Define $\Phi_\phi (t):=\|z_t(\phi)\|$. Since (5.2) is an autonomous system, then $\Phi_\phi (t_2)=\Phi_{z_{t_1}(\phi )} (t_2-t_1)$
for $t_2>t_1>0$, and the above estimate proves that
  $\Phi_\phi (t_2)<  \Phi_\phi (t_1)$ if $t_2> t_1+\tau$. The same arguments yield that $t\mapsto \|z_t(\phi)\|$ is non-increasing for $t\ge 0$, so $\Phi_\phi(t)\searrow \al$ as $t\to\infty$, for some $\al\ge 0$.

 Next, 
 consider the $\omega$-limit set $\omega(\phi)$, which is non-empty. The invariance of
 $\omega(\phi)$ under (5.2) implies that $\omega(\phi)\subset \{\psi \in \bar S: \| \psi\| =\al\},$ where $\bar S$ denotes the closure of $S$ in $C$. But the components $z_i(t)$ are bounded away from $-1$ (cf. Theorem 2.2), and therefore $\omega(\phi)\subset S$.
 
  
 If $\al>0$, let $\psi\in \omega(\phi)$. We have $\psi\in S$ and $\|\psi\|=\al$.
  However this is not possible, since $z_t(\psi)\in\omega(\phi)$ and from (5.3) we get $\|z_t(\psi)\|<\|\psi \|=\al$ for $t\ge \tau.$
 This shows that $\al=0$, and the theorem is proved.\ter

 \med
 
 {\bf Remark 5.1}. In [3], the  global asymptotic stability (with respect to $C^+_0$) of $x^*$ was proved under the stronger hypothesis
 $1<\gamma_i\le \min \{e^2,e^{x_i^*}\},\,  i=1,\dots,n ,$ which turned out to be very restrictive, since for $x_i^*=\max_{1\le j\le n} x_j^*$ we necessarily have  $\gamma_i\ge e^{x_i^*}$, and where the equality holds if and only if either $a_{ij}=0$ or $x_j^*=x_i^*$ for all $j\ne i$. Furthermore, criteria for the existence of such a positive equilibrium were not established in [3].

 \med
 
 In the above proof, observe that hypothesis {\bf (A2)} was not directly applied to system (5.2), obtained as a consequence of  the change of variables (5.1). Actually, {\bf (A2)} was used  only to guarantee the existence of a positive equilibrium with all its components in the interval $(0,2]$, which is crucial for two reasons:  on one hand, its local stability is deduced regardless of the size of the positive delays, and, on the other hand, Lemma 5.1 can be applied. Note that the estimate in Lemma 5.1 is no longer valid for $x>2$. This observation permits to state the global attractivity of the positive equilibrium under weaker assumptions, as follows.
 
 \proclaim {Theorem 5.2}. Assume (2.6)
 for some positive vector $c=(c_1,\dots,c_n)$. Then, the unique positive equilibrium $x^*=(x_1^*,\dots, x_n^*)$ (whose existence is given by Theorem 4.2) 
  is GAS if $x_i^*\le 2$ for $i=1,\dots,2$.

 {\bf Remark 5.2}. For the scalar Nicholson's blowflies equation, it is well-known that if $\gamma_1=\be_1/d_1>e^2$,  large delays can lead to the existence of  periodic solutions appearing from a  Hopf bifurcation. Also for $n>1$, we can show that hypothesis {\bf (A2)} is  a sharp condition for the {\it absolute}  global asymptotic stability  (i.e., for the global asymptotic stability independently of the size of positive delays $\tau_{ik}$) of $x^*$;   if $\gamma_i>e^2$ for some $i$, in general large delays bring instability, as illustrated in the example below.

 \med

 {\bf Example 5.1}. Consider (1.1) with $n=2$, $m=1$:
  $$
 \eqalign{
x_1'(t)&=-d_1x_1(t)+a_{12}x_2(t)+\be_1x_1(t-\tau_1)e^{-x_1(t-\tau_1)}\cr
x_2'(t)&=-d_2x_2(t)+a_{21}x_1(t)+\be_2x_2(t-\tau_2)e^{-x_2(t-\tau_2)}\cr}\eq(5.5)
$$
and $a_{12}\ge 0, a_{21}, d_i,\be_i,\tau_i>0, i=1,2,$ with $1<\gamma_1=\be_1/(d_1-a_{12})\le e^2$ and $ \gamma_2=\be_2/(d_2-a_{21})>e^2$, so that {(\bf A2)} fails.
  Under some further conditions on the coefficients in (5.5),  we   show that the positive equilibrium $x^*=(x_1^*,x_2^*)$ is not asymptotically stable if the size of the delay $\tau_2$ is large.

Let   $a_{21}>0$ be sufficiently small  so that $\be_2/d_2>e^2$. The linearization about  $x^*=(x_1^*, x_2^*)$ is given by 
$$y_i(t)=-[d_iy_i(t)+L_{i1}(y_t)+L_{i2}(y_t)],\q i=1,2,$$
where the linear operators $L_{ij}$ are defined by
$$\eqalign{
L_{11}(\var)&=-\be_1h'(x_1^*)\var_1(-\tau_1),\ L_{12}(\var)=0\cr
L_{21}(\var)&=-a_{21}\var_1(0),\  L_{22}(\var)=-\be_2h'(x_2^*)\var_2(-\tau_2),\q \var=(\var_1,\var_2)\in C.\cr}$$
Define now 
$$\hat N=D-\Big [\| L_{ij}\|\Big]=\pmatrix{d_1-\be_1|h'(x_1^*)|&0\cr -a_{21}&d_2-\be_2|h'(x_2^*)|\cr},$$
with eigenvalues $\la_1=d_1-\be_1|h'(x_1^*)|$ and $\la_2=d_2-\be_2|h'(x_2^*)|$. We  claim that it is possible to have $\la_2<0$. If this is the case,  from Theorem 2.3 in [4] we conclude that there is $\tau_2>0$ for which the equilibrium $x^*=(x_1^*,x_2^*)$ of (5.5) is unstable.

For  $\al:=a_{21}x_1^*$, we have $e^{x_2^*}={{\be x_2^*}\over {d_2 x_2^*-\al}}\to \be_2/d_2>e^2$ as $\al\to 0^+$. This implies $x_2^*=x_2^*(\al)>2$, for either $a_{21}$ or $x_1^*$ small (for instance, with $a_{12}=0$, we have that $x_1^*=\log (\be_1/d_1)\to 0^+$ if $\be_1/d_1\to 1^+$). Thus, $\la_2=d_2+\be_2(1-x_2^*)e^{-x_2^*}$ and for $x_2^*(0):=\log (\be_2/d_2)$
we obtain
$$\la_2=\la_2(\al)={1\over {x_2^*}}[-d_2({x_2^*})^2+(2d_2+\al){x_2^*}-\al]\to d_2(2-x_2^*(0))<0\q {\rm as}\q \al\to 0^+.$$

{
\input epsf
\centerline{\epsfxsize=6cm
(a) \epsfbox{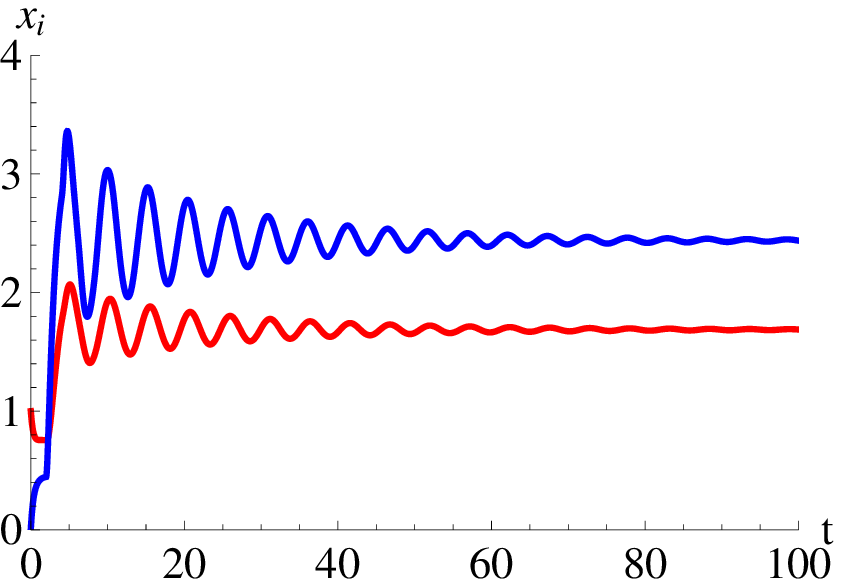}
\epsfxsize=6cm
\hskip0.2in (b) \epsfbox{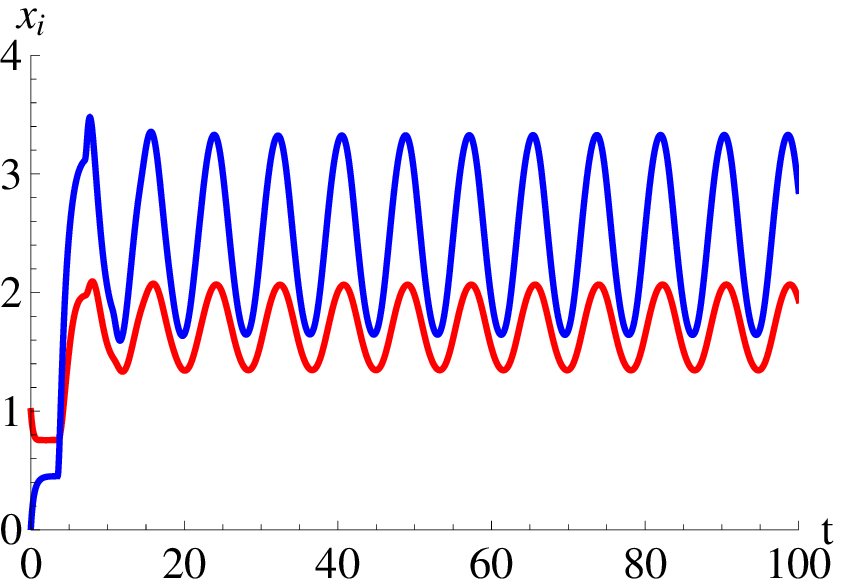}}
\noindent {Figure 3. Illustration of Example 5.1. Parameters are $a_{12}=a_{21}=1$, $d_1=d_2=2$, $\beta_1=3$, $\beta_2=15$, $\tau_1=1$. Then 
$\gamma_1=\beta_1<e^2$ and $\gamma_2=\beta_2>e^2$. In (a), we set $\tau_2=2$,
and we observe the convergence of solutions to an equilibrium. Increasing the delay to $\tau_2=3.5$, the equilibrium becomes unstable and we can see a periodic oscillation in (b).}
}

\

{\bf Acknowledgements}: 
Work supported by Funda\c c\~ao para a Ci\^encia e a Tecnologia,  PEst-OE/MAT/UI0209/2011 (T. Faria) and by ERC Starting Grant Nr. 259559
and ESF project FuturICT.hu 
(T\' AMOP-4.2.2.C-11/1/KONV-2012-0013) (G. R\"ost).

\

{\bf \centerline{References}}

\med

\baselineskip=13.5pt

\item{1.}  Berezansky, L., Idels,  L.,  and  Troib, L.  (2011). Global dynamics of Nicholson-type delay systems with applications.
Nonlinear Anal. Real World Appl. 12, 436--445.

\item{2.}   Berezansky, L., Braverman, E., Idels, L. (2010). Nicholson's blowflies differential equations revisited: main results and open problems. Appl. Math. Model. 34  1405--1417.

\item {3.}  Faria, T. (2011). Global asymptotic behaviour for a Nicholson model
 with patch structure and multiple delays. Nonlinear Anal. 74, 7033--7046.
 
 \item{4.} Faria, T., and  Oliveira, J. J. (2008).
Local and global stability for   Lotka-Volterra systems 
with distributed delays and instantaneous feedbacks. 
 J. Differential Equations 244, 1049--1079.

 \item{5.} Fiedler, M.  (1986). Special Matrices and Their Applications in Numerical Mathematics, Martinus Nijhoff
Publ. (Kluwer), Dordrecht.

\item {6.}  Gurney, W. S. C.,  Blythe S. P., and  Nisbet, R. M. (1980). Nicholson's blowflies revisited. Nature  287, 17--21.

\item{7.}  Hale, J. K. (1988). Asymptotic Behavior of Dissipative Systems,
Amer.~Math.~Soc., Providence, Rhode Island.

\item{8.}  Hofbauer, J. (1990). An index theorem for dissipative systems. Rocky Mountain J. Math.  20, 1017--1031.

\item{9.} Kuang,  Y. (1993). Delay Differential Equations with Applications in Population Dynamics, Academic
Press, London.

\item{10.}  Liu, B. (2009). Global stability of a class of delay differential systems. J. Comput. Appl. Math. 233, 217--223.

\item{11.}  Liu, B. (2010). Global stability of a class of Nicholson's blowflies model with patch structure and multiple time-varying delays. Nonlinear Anal. Real World Appl. 11, 2557--2562.

\item {12.}  Liu, X., and Meng, J. (2012). The positive almost periodic solution for Nicholson-type delay systems with linear harvesting term. Appl. Math. Model. 36,  3289--3298.

\item {13.}  Nicholson, A. J. (1954). An outline of the dynamics of animal populations. Austral. J. Zool. 2, 9--65.

\item {14.}  R\"ost, G., and  Wu, J. (2007). Domain-decomposition method for the global dynamics of delay differential equations with unimodal feedback. Proc.~R.~Soc.~Lond.~Ser.~A ~Math.~Phys.~Eng. Sci. 463, 2655--2669

\item{15.}  Smith, H. L.  (1995). Monotone Dynamical Systems. An
Introduction to the Theory of Competitive and Cooperative Systems,
 Mathematical Surveys and Monographs, Amer.~Math.~Soc., Providence, RI.
 
  \item{16.}  Smith, H. L., and  Thieme, H. R. (2011). Dynamical Systems and Population Persistence,  
Amer.~Math.~Soc., Providence, RI.

 \item{17.}  Smith, H. L., and Waltman, P. (1995).  The Theory of the
Chemostat,  University Press, Cambridge.


 \item{18.}  Wang, L. (2013). Almost periodic solution for Nicholson's blowflies model with patch structure and linear harvesting terms. Appl. Math. Model. 37,  2153--2165.
 

\item {19.}  Zhao, X.-Q., and  Jing, Z.-J. (1996). Global asymptotic behavior in some cooperative systems of functional differential equations.
Cann. Appl. Math. Quart. 4 , 421--444.

\end